\documentclass{article} 
\usepackage[cp1251]{inputenc}
\usepackage{latexsym,amsfonts,amssymb,amsmath}
\usepackage{euscript}
\usepackage[T2A]{fontenc}
\usepackage{amsthm}
\usepackage[english]{babel}
\usepackage{mathrsfs}

mm \oddsidemargin=23mm \evensidemargin=23mm \headsep=6.5mm
\begin{document}


\noindent UDC\  519.41/47

\noindent {\bf O.O.~Trebenko}

\noindent{\footnotesize (Institute of Mathematics of the National
Academy of Sciences of Ukraine, Kyiv)}

\vspace{2mm}\noindent {\bf\large ON GROUPS WITH A SUPERCOMPLEMENTED
SUBGROUP}

\vskip 2mm

\noindent{\small Groups, in which every subgroup containing some
fixed primary cyclic subgroup has a complement, are investigated}

\medskip

\setcounter{page}{1}

\newcommand{\ang}[1]{\langle{#1}\rangle}
\selectlanguage{english}

\thispagestyle{empty}

{\bf 1. Introduction.} 
Recall that a subgroup $H$ of the group $G$ is called complemented
in $G$, if there exists some subgroup $T$ of $G$ such that $G=HT$
and $H\cap T=1$; it is said that $T$ complements $H$ in $G$ and $T$
is a complement to $H$ in $G$. Finite groups in which all subgroups
are complemented were first considered by Ph.Hall \cite{Hall}. A
complete constructive description of arbitrary groups in which all
subgroups are complemented has been obtained by N.V.Chernikova
\cite{NVCher,NCher}. In \cite{NVCher} such groups were called
completely factorizable. In view of N.V.Chernikova's Theorem
\cite{NVCher,NCher}, completely factorizable groups are solvable
(more precisely, metabelian) and locally finite.

Following \cite{Kreknin} we call a subgroup $H$ of the group $G$
supercomplemented in $G$, if each subgroup of $G$ containing $H$ is
complemented in $G$. In connection with N.V.Chernikova's Theorem, it
is naturally to investigate groups which have a proper
supercomplemented subgroup. In \cite{Ch_Tr1} locally graded
$p$-groups with a supercomplemented cyclic subgroup were
investigated. According to \cite{Ch_Tr1}, such groups are locally
finite and solvable, and also, in the case when $p\neq 2$,
metabelian. Later, V.A.Kreknin \cite{Kreknin} has showed that the
derived length of such 2-groups does not exceed 3.

Recall that a group is said to be locally graded, if its every
nonidentity finitely generated subgroup has a proper subgroup of
finite index \cite{lg}. The class of all locally graded groups is
extremely wide. For example, all locally finite, solvable and
locally solvable, residually finite groups, linear groups, radical
(in the sense of B.I.Plotkin) groups, $RN$-groups (and, at the same
time,  groups of all Kurosh-S.N.Chernikov's classes) are locally
graded.

(Note that by a proper subgroup of the group $G$ we mean, as in
\cite{Ch_Tr2}, a subgroup different from $G$).

\vskip2mm The main results of the present paper are the following
Theorems 1--4.

Below, as usual, for real number $r\geq 0$, $[r]$ is the greatest
integer not exceeding $r$.

\vskip2mm \textbf{Theorem 1.} \textit{Let $G$ be an $RN$-group with
a supercomplemented cyclic $p$-subgroup $\ang{x}$ of order $m$.
Then: \vspace{-1mm}
\begin{itemize}\itemsep-0.1cm
  \item [(i)] $G$ is solvable and locally finite.
  \item [(ii)] In the cases when $m=1$, $m=2$, $2< m<8$ and $m\geq 8$,
  for the derived length $d(G)$ of $G$, respectively, $d(G)\leq 2$,
  $d(G)\leq 11$, $d(G)\leq 18$ and $d(G)\leq 5\log_{9}\frac{n-2}{8}+13$ ($> 18$),
  where
$n=[m(m-1)+m\log_{2}m]$.
  \item [(iii)] $G$ is residually finite.
  \item [(iv)] Each $p$-subgroup of
  $G$ is nilpotent, almost elementary abelian and solvable
  of derived length $\leq 3$. In the case when $p\neq 2$, each
  $p$-subgroup of $G$ is metabelian.
\end{itemize}}

\textbf{Theorem 2.} \textit{Let $G$ be a periodic locally graded
group with a supercomplemented cyclic $p$-subgroup of order $m$ and
also $2\notin\pi(G)$. Then the statements (i)--(iv) from Theorem 1
are valid.}

\vskip2mm It is naturally to consider groups in which all subgroups
not contained in some proper fixed subgroup are complemented. As
usually, we call such proper subgroup a $C$-separating subgroup. In
\cite{Kr_Sp} the theorem on solvability of a finite group with a
$C$-separating subgroup was established. It is established
\cite{Ch_Tr2} that a locally almost solvable group with a
$C$-separating subgroup is locally (solvable and finite). Remark
that in contrast to \cite{Kr_Sp}, the fundamental theorems of the
Theory of finite simple groups were not used in \cite{Ch_Tr2}.

The following theorem, in particular, considerably generalizes
Theorem \cite{Ch_Tr2} mentioned above.

\vskip2mm \textbf{Theorem 3.} \textit{Let $G$ be a locally graded
group with a $C$-separating subgroup $H$. Then:
\vspace{-1mm}\begin{itemize}\itemsep-0.1cm
  \item [(i)] $G$ is solvable and locally finite.
  \item [(ii)] $G$ is residually finite.
  \item [(iii)] For some $p\in \mathbb{P}$, $G$ possesses a cyclic
  supercomplemented $p$-subgroup, each $p$-subgroup of
  $G$ is nilpotent, almost elementary abelian and solvable of derived length $\leq 3$,
   and each $q$-subgroup of $G$ with $q\neq p$
  is elementary abelian.
\end{itemize}
}

The proof of Theorem 3 essentially uses Theorem 1.

\vskip2mm  Finally, Theorem 4 shows that the class of all groups
containing a $C$-separating subgroup is a proper subclass of the
class of all groups containing a supercomplemented primary cyclic
subgroup.

\vskip2mm \textbf{Theorem 4 (N.S.Chernikov, O.O.Trebenko).}
\textit{The class $\mathfrak{X}$ of all groups containing a
$C$-separating subgroup is a proper subclass of the class
$\mathfrak{Y}$ of all groups containing a supercomplemented primary
cyclic subgroup. (More concretely, $\mathfrak{X}$ is contained in
$\mathfrak{Y}$ and the holomorph $G$ of the cyclic group of order 8
belongs to $\mathfrak{Y}\backslash \mathfrak{X}$). }

\vskip2mm In what follows, $\mathbb{N}$ and $\mathbb{P}$ denote the
sets of all naturals and primes, respectively. The symbol
$\leftthreetimes$ is used to denote the semidirect product. Let $G$
be a group. Notations $H \leq G$ ($H\trianglelefteq G$) and $H<G$
($H\vartriangleleft G$) mean that $H$ is a (normal) subgroup of $G$
and $H$ is a (normal) subgroup of $G$ different from $G$,
respectively. $\pi(G)$ is the set of all $p\in \mathbb{P}$ for which
$G$ has an element of order $p$. Further, $\Phi(G)$ is the Frattini
subgroup of $G$, $Z(G)$ is the centre of $G$. Let $\varnothing\neq
H\subseteq G$. Then $H^{G}=\{h^{g}:h\in H, \:g\in G\}$,
$H^{n}=\ang{h^{n}:h\in H, \:n\in \mathbb{N}}$. Below, in the case
when $G$ is solvable, $d(G)$ is its derived length. Other notations
are standard.

\vskip2mm In connection with Theorem 1, remark that for $p\neq 2$, a
locally graded $p$-group with a supercomplemented cyclic subgroup
$\ang{x}$ is metabelian, nilpotent and $G=\ang{x}B$, $\ang{x}\cap
B=1$, where $B$ is elementary abelian (see Theorem 2 \cite{Ch_Tr1}).
It may be that $\ang{x},B\ntrianglelefteq G$, as shows the following
example.

\vskip2mm \textbf{Example.} Let $G=A\leftthreetimes F$ be a finite
$p$-group, where $A=\ang{b}\times\ang{c}$, $|\ang{b}|=|\ang{c}|=p$,
$F=\ang{x}\leftthreetimes \ang{a}$, $|\ang{x}|=p^{2}$,
$|\ang{a}|=p$, and $x^{a}=x^{p+1}$  and also $b^{x}=bc$, $c^{x}=c$,
$b^{a}=b$, $c^{a}=c$. Then $G$ may be presented in the form:
$G=\ang{x}B$, $\ang{x}\cap B=1$, where
$B=\ang{a}\times\ang{b}\times\ang{c}$. Obviously,
$\ang{x}\ntrianglelefteq G$ and $B\ntrianglelefteq G$.

\vskip2mm Main ideas, methods and approaches used in the present
paper are worked up in the paper \cite{Ch_Tr3} submitted to
publishing earlier.

\vskip3mm 
{\bf 2. Preliminary results.}

\vskip2mm\textbf{Lemma 1.} \textit{Let $G$ be a group with a
supercomplemented subgroup $H$ and $H\subseteq K\leq G$,  and let
$\varphi$ be a homomorphism of $K$. Then $H^{\varphi}$ is
supercomplemented in $K^{\varphi}$.}

\vskip2mm {\bfseries \itshape Proof.} Take any $S\leq K^{\varphi}$
such that $H^{\varphi}\subseteq S$. Let $L\leq K$ and
$\mathrm{Ker}\, \varphi\subseteq L$, $L^{\varphi}=S$. Then
$H\subseteq L$. Take a complement $T$ to $L$ in $K$. Then $K
^{\varphi}=S T ^{\varphi}$ and, obviously, $S \cap T ^{\varphi}=1$.
So $T ^{\varphi}$ is a complement to $S$ in $K ^{\varphi}$. Lemma is
proven.

\vskip2mm\textbf{Proposition 1.} \textit{Let $G$ be a group with a
supercomplemented cyclic $p$-subgroup $\ang{x}$ of order $m$. If for
some $q\in \mathbb{P}$, an elementary abelian $q$-subgroup $Q$ of
$G$ is its minimal normal subgroup, then either $m\neq 1$ and
$|Q|\leq q^{(m-1)m}m^{m}$, or $m=1$ and $|Q|= q$.}

\vskip2mm{\bfseries \itshape Proof.} Let $m\neq 1$. Show that $Q$
has a minimal normal subgroup of order $\leq q^{m-1}$ of the group
$Q\ang{x}$. Take $a\in Q\backslash\{1\}$. Denote $K=\ang{a^{x},
a^{x^{2}},\ldots, a^{x^{m}}}$. Clearly, $K\trianglelefteq
{Q\ang{x}}$ and $|K|\leq q^{m}$. Obviously, $K$ contains some
minimal normal subgroup $Q_{1}$ of $Q\ang{x}$.

Suppose that $|Q_{1}|= q^{m}$. Then
$Q_{1}=K=\ang{a^{x}}\times\ang{a^{x^{2}}}\times\ldots\times
\ang{a^{x^{m}}}$. Clearly, $\ang{a^{x}a^{x^{2}}\ldots a^{x^{m}}}\lhd
Q\ang{x}$ and $|\ang{a^{x}a^{x^{2}}\ldots a^{x^{m}}}|=q<q^{m}$,
which is a contradiction. Thus, $|Q_{1}|\leq q^{m-1}$.

Further, for some $T\leq G$, $G=(Q_{1}\ang{x})T$ and $(Q_{1}
\ang{x})\cap T=1$. In view of S.N.Chernikov's Lemma (see, for
instance, \cite{SNLemma}, Lemma 1.8),
$Q\ang{x}=Q_{1}\ang{x}(Q\ang{x}\cap T)$. We have
$Q_{1}\ang{x}\cap(Q\ang{x}\cap T)=1$. Denote $Q_{2}=Q\cap T$. Then
$Q_{2}\neq Q$. Since, obviously, $G=Q\ang{x}T=\ang{x}QT=QT\ang{x}$
and $Q_{2}\trianglelefteq Q,T$, we have $G=N_{G}(Q_{2})\ang{x}$. So
for each $g\in G$, there exists $h\in \ang{x}$ such that
$Q_{2}^{g}=Q_{2}^{h}$. Further, for any $h\in \ang{x}$, obviously,
\begin{gather*}
    |Q:Q_{2}^{h}|=|Q:Q_{2}|\leq|Q\ang{x}:Q\ang{x}\cap T|=|Q_{1}\ang{x}|\leq q^{m-1}m.
\end{gather*}
So \[|Q:{\mathop {\cap} \limits_{h\in \ang{x}}} Q_{2}^{h}|\leq
{\mathop {\prod} \limits_{h\in \ang{x}}} |Q:Q_{2}^{h}|\leq
(q^{m-1}m)^{m}.\] As far as ${\mathop {\cap} \limits_{h\in \ang{x}}}
Q_{2}^{h}={\mathop {\cap} \limits_{g\in G}} Q_{2}^{g}\trianglelefteq
G$, $Q_{2}\neq Q$ and Q is a minimal normal subgroup of $G$,
${\mathop {\cap} \limits_{h\in \ang{x}}} Q_{2}^{h}=1$ and so
$|Q|\leq (q^{m-1}m)^{m}$.

If $m=1$, then $G$ is completely factorizable. Therefore, in
consequence of N.V.Chernikova's Theorem \cite{NVCher,NCher},
$|Q|=q$. Proposition is proven.

\vskip2mm\textbf{Remark.} In view of Zassenhaus'es Theorem, for an
arbitrary $n\in \mathbb{N}$, the derived lengths of solvable linear
groups of degree $\leq n$ over fields are bounded by some natural
number depending only on $n$ (see, for instance, \cite{Zass},
Theorem 3.7). Let $\zeta(n)$ be the smallest such number. Obviously,
$\zeta(n)\leq \zeta(n+1)$.

In 1958 B.Huppert (see, for instance, \cite{Merzlyakov}, Theorem
45.2.1) has showed that $\zeta(n)\leq 2n$. For $n\geq 66$, M.Newman
\cite{Newman} has obtained the following estimation: $\zeta(n)\leq
5\log_{9}\frac{n-2}{8}+10$.

Consequently, for $n\leq 6$, we have $\zeta(n)(\leq 2n)\leq 12$, for
$7\leq n\leq 73$, $\zeta(n)\leq 14$, and for $n> 73$, we have
$\zeta(n)\leq [5\log_{9}\frac{n-2}{8}+10](\geq 15)$.

\vskip3mm\textbf{Proposition 2.} \textit{Let $G$ be a finite
solvable group with a supercomplemented cyclic $p$-subgroup
$\ang{x}$ of order $m$, and let $n=[m(m-1)+m\log_{2}m]$ if $m\neq
1$. Then $d(G)\leq \zeta(n)+3$. Moreover, in the cases when $m=1$,
$m=2$, $2<m<8$ and $m\geq 8$ the following holds respectively:
$d(G)\leq 2$, $d(G)\leq 11$, $d(G)\leq 18$ and $d(G)\leq
5\log_{9}\frac{n-2}{8}+13(> 18)$.}

\vskip2mm{\bfseries \itshape Proof.} We may assume, of course, that
$G\neq 1$. Let $G=G_{0}\supset G_{1}\supset \ldots\supset G_{s}=1$
be some chief series of $G$. For some $q=q(i)\in \mathbb{P}$,
$G_{i}/G_{i+1}$ is an elementary abelian $q$-group. Since
$G_{i}/G_{i+1}$ is a minimal normal subgroup of $G/G_{i+1}$  and
$\ang{x}G_{i+1}/G_{i+1}$ is a supercomplemented cyclic $p$-subgroup
of order $\leq m$ of $G/G_{i+1}$ (see Lemma 1), $|G_{i}/G_{i+1}|\leq
q^{n}$ where $n=1$ if $m=1$, and $n$ is as above if $m\neq 1$. (see
Proposition 1).

Then, since $G/C_{G}(G_{i}/G_{i+1})$ is isomorphically embedded into
$GL_{n}(q)$, by Zassenhaus'es Theorem,
$d(G/C_{G}(G_{i}/G_{i+1}))\leq \zeta(n)$ (see Remark above).
Therefore,\quad $d(G/{\mathop {\bigcap} \limits_{i=0}^{s-1}}
C_{G}(G_{i}/G_{i+1}))\leq \zeta(n)$.

Let $F={\mathop {\bigcap} \limits_{i=0}^{s-1}}
C_{G}(G_{i}/G_{i+1})$. In view of Theorem III.4.3 \cite{Huppert}
(for instance), $F$ is nilpotent.

Let $F\neq 1$, $q\in \pi(F)$ and $P$ be the Sylow $q$-subgroup of
$F$. By Lemma 1, $\ang{x}$ is a supercomplemented subgroup of the
group $P\ang{x}$. Therefore, in view of Theorem 2 \cite{Ch_Tr1} and
Theorem \cite{Kreknin}, $d(P\ang{x})\leq 3$ if $q=p$.

Let $q\neq p$. For some $T\leq Q\ang{x}$, $Q
\ang{x}=(\Phi(Q)\ang{x})T$ and $(\Phi(Q)\ang{x})\cap T=1$.
Obviously, $T$ is a complement to $\Phi(Q)$ in $Q$. Consequently,
$\Phi(Q)=1$ and $Q$ is elementary abelian. So $d(Q)=1$.

Since $F$ is decomposed into the direct product of its primary Sylow
subgroups, $d(F)\leq 3$. Therefore, $d(G)\leq \zeta(n)+3$.

Further, it is easy to see that $n<66$ iff $m\leq 7$.

If $m=1$, then $G$ is completely factorizable. Therefore, in view of
N.V.Chernikova's Theorem, $d(G)\leq 2$.

If $m=2$, then $n=4$ and, with regard to Remark, $d(G)\leq 2n+3=11$.

If $m=8$, then $n=80$. So, with regard to Remark, for $m\geq 8$,
$d(G)\leq 5\log_{9}\frac{n-2}{8}+13$, and for $m=8$, $d(G)\leq
\zeta(80)+3\leq [5\log_{9}\frac{80-2}{8}+13]=18$.

If $2<m<8$, then, with regard to Remark, $d(G)\leq \zeta(n)+3\leq
\zeta(80)+3\leq 18$. Proposition is proven.

\vskip4mm\textbf{Proposition 3.} \textit{Let $G$ be a locally
solvable group with a finite supercomplemented subgroup $\ang{x}$.
Then $G$ is locally finite.}

\vskip2mm{\bfseries \itshape Proof.}  Obviously, it is sufficient to
show that every finitely generated subgroup of $G$ containing
$\ang{x}$ is finite. Thus, with regard to Lemma 1, the proof is
reduced to the case when $G$ is finitely generated solvable. Since
$G$ is solvable, it has a finite normal series with abelian factors.
Inasmuch as every abelian group is periodic-by-torsion-free, this
series is contained in some finite normal series $\mathcal{M}$ of
$G$ with periodic abelian and torsion-free abelian factors.

Assume that $G$ is infinite. Then, in view of S.N.Chernikov's
Theorem (see, for instance, \cite{SNProp}, Proposition 1.1), $G$ is
non-periodic. Therefore, for some neighbouring terms $N$ and
$L\subset N$ of series $\mathcal{M}$, $N/L$ is torsion-free and
$G/N$ is periodic. In consequence of Proposition 1.1 \cite{SNProp},
$G/N$ is finite. Further, by virtue of Lemma 1, $\ang{x}L/L$ is
finite supercomplemented subgroup of $G/L$. Thus, the proof of the
present theorem is reduced to the case when for some torsion-free
abelian subgroup $N\unlhd G$, $|G:N|<\infty$.

By Schreier-Dyck's Theorem (see, for instance, \cite{Kurosh},
pp.228,111), $N$ is finitely generated. Therefore, for any $q\in
\mathbb{P}$, $N^{q}\neq N$.

Clearly, $|\pi(G)|<\infty$. Take $q\in\mathbb{P}\setminus\pi(G)$.
According to Lemma 1, $\ang{x}N^{q}$ is complemented in $\ang{x}N$
by some subgroup $H$. It is easy to see that $|H|=|N:N^{q}|=q^{n}$
where $n\in \mathbb{N}$. So $H$ is a nonidentity $q$-subgroup of
$G$, which is a contradiction. Proposition is proven.

\vspace{2mm}Recall that an involution of a group is its element of
order 2.

\vskip2mm\textbf{Proposition 4.} \textit{The periodic group $G$
without involutions is locally graded iff it is an $RN$-group.}

\vskip2mm{\bfseries \itshape Proof.}  Let $G$ be a periodic group
without involutions. If $G$ is an $RN$-group, then it is locally
graded. Suppose that $G$ is locally graded and $G\neq 1$. Let $K$ be
any nonidentity finitely generated subgroup of $G$. Then, with
regard to Poincare's Theorem, there exists $N\vartriangleleft K$
with $|N:K|<\infty$. Since $K/N$ is a nonidentity finite group and,
obviously,  $2\notin\pi(K/N)$, in view of Feit-Thompson's Theorem
\cite{Feit}, $K/N$ is solvable. Hence $(K/N)'\neq K/N$. Therefore
$K'\neq K$. So, by S.Brodski\v{\i}'s Theorem (see, for instance,
\cite{Brodskii}), $G$ is an $RN$-group. Proposition is proven.

\vskip2mm\textbf{Lemma 2.} \textit{Let for some $m\in \mathbb{N}$, a
$p$-group $G$ has some local system of subgroups
$\{G_{\iota}:\iota\in I\}$ such that each $G_{\iota}$ has a (normal)
elementary abelian subgroup of index $\leq m$. Then $G$ has a
(normal) elementary abelian subgroup of index $\leq m$ and $G$ is
nilpotent. In particular, $G$ is almost elementary abelian.}

\vskip2mm{\bfseries \itshape Proof.}  By Poincare's Theorem,
$G_{\iota}$, obviously, has a normal elementary abelian subgroup of
finite index and, at the same time, locally finite subgroup of
finite index. Therefore, in view of O.J.Schmidt's Theorem (see, for
instance, \cite{Kurosh}, p.337), $G$ is locally finite. Then each
$G_{\iota}$ and, at the same time, $G$ have a local system of finite
subgroups with a (normal) elementary abelian subgroup of finite
index $\leq m$. Taking this into account, further we may assume
without loss of generality that all $G_{\iota}$ are finite.

Let $M_{\iota}$ be a set of all (normal) elementary abelian
subgroups of index $\leq m$ of $G_{\iota}$. Let $M_{\alpha} \leq
M_{\iota}$ iff $G_{\alpha} \subseteq G_{\iota}$. In the case
$M_{\alpha} \leq M_{\iota}$ we define the projection $\pi_{\iota
\alpha}$ from $M_{\iota}$ into $M_{\alpha}$ as follows: for an
arbitrary $K \subseteq M_{\iota}$ $K^{\pi_{\iota \alpha}}=K \cap
G_{\alpha}$. Obviously, the following holds:

\vspace{-1mm}\begin{enumerate}\itemsep-0.09cm
    \item[1)] for each $M_{\alpha}$ and $M_{\beta}$, there exists
    $M_{\gamma}$ such that $M_{\alpha},M_{\beta} \leq M_{\gamma}$;
    \item[2)] if $M_{\alpha} \leq M_{\beta}$, $M_{\beta} \leq
    M_{\gamma}$, then $\pi_{\gamma \alpha}=\pi_{\gamma \beta}\pi_{\beta
    \alpha}$;
    \item[3)] $\pi_{\iota \iota}$ is an identity mapping of
    $M_{\iota}$ onto itself.
\end{enumerate}

\vspace{-1mm}Consequently, in view of (\cite{Kurosh}, p.351-353),
there exist $K_{\iota}\in M_{\iota}$, $\iota\in I$, such that
$K_{\alpha}=K_{\iota}\cap G_{\alpha}$ whenever $G_{\alpha} \subseteq
G_{\iota}$. Obviously, $K=\underset{\iota \in I}{\bigcup} K_{\iota}$
is elementary abelian.

Let $n=\underset{\iota \in I}\max|G_{\iota}:K_{\iota}|$ and let for
$\gamma \in I$, $|G_{\gamma}:K_{\gamma}|=n$ and
$G_{\gamma}=\overset{n}{\underset{i=1}{\bigcup}} a_{i}K_{\gamma}$.
Take any $g\in G_{\iota}$. Then for some $\beta\in I$, $g\in
G_{\beta}\supseteq G_{\gamma}$. Since $K_{\gamma}\subseteq
K_{\beta}$ and $|G_{\beta}:K_{\beta}|\leq |G_{\gamma}:K_{\gamma}|$,
obviously, $G_{\beta}=\overset{n}{\underset{i=1}{\bigcup}}
a_{i}K_{\beta}$. Consequently, $g\in \overset{n}{\underset
{i=1}{\bigcup}} a_{i}K$. Thus, $G=\overset{n}
{\underset{i=1}{\bigcup}} a_{i}K$. At the same time, $|G:K|\leq
n\leq m$. In consequence of Poincare's Theorem, $|G:{\mathop {\cap}
\limits_{g\in G}} K^{g}|<\infty$. So $G$ is almost elementary
abelian.

Suppose that each $M_{\iota}$ consists of all normal elementary
abelian subgroups of $G_{\iota}$. Take any $g\in G$ and $a\in K$.
Then for some $\alpha\in I$, $g\in G_{\alpha}$ and $a \in
K_{\alpha}\unlhd G_{\alpha}$. Therefore $a^{g} \in
K_{\alpha}\subseteq K$. Thus $K\unlhd G$.

It is clear that an exponent of $G$ is finite. Therefore, in view of
Baumslag's Theorem \cite{Baumslag}, $G$ is nilpotent. Lemma is
proven.

\vspace{3mm}
{\bf 3. Proofs of theorems.}

\vspace{2mm} {\bfseries\itshape Proof of Theorem 1}. (i) It is
sufficient to show, with regard to Lemma 1, that the theorem is
valid for any finitely generated subgroup $G^{*}\supseteq \ang{x}$
of $G$. Thus the proof is reduced to the case when $G$ is finitely
generated.

Let $K$ be an intersection of all $N\trianglelefteq G$ such that
$G/N$ is finite and solvable. By virtue of Proposition 2 and Lemma
1, $d(G/N)\leq \zeta(n)+3$, where $n=\max\{[m(m-1)+m\log_{2}m],
1\}$. Consequently, $G/K$ is solvable with $d(G/K)\leq \zeta(n)+3$.
Further, in view of Proposition 3 and Lemma 1, $G/K$ is finite.

Suppose $K\neq 1$. By Schreier-Dyck's Theorem, an $RN$-group $K$ is
finitely generated. Let $K=\ang{M}$ where $|M|<\infty$,
$\mathcal{A}$ be a series with abelian factors of $K$ and $U$ be a
union of all $T\in \mathcal{A}\backslash \{K\}$. Since $M\nsubseteq
T$ and $M$ is finite, clearly, $M\nsubseteq U$. So $U\neq K$.
Clearly, $U\vartriangleleft K$ and $K/U$ is abelian. So $K\neq K'$.

If $|K:K'|<\infty$, then $G/K'$ is finite and solvable. Therefore
$K'\supseteq K$, which is a contradiction. So $|K:K'|$ is infinite.

Take any $q\in \mathbb{P}$. Since $K/K'$ is infinite abelian
finitely generated, $(K/K')^{q}\neq K/K'$. Let $(K/K')^{q}=L/K'$.
Then $L\lhd G$, $L\subset K$ and $K/L$ is finite abelian. So $G/L$
is finite and solvable, which is a contradiction. Thus $K=1$.
Therefore $G$ is finite and, with regard to Proposition 2, (i) is
valid.

(ii) follows from (i) and Proposition 2 (with regard to Lemma 1).

(iii) Let $1\neq g\in G$. Consider $R=\ang{x,g}$. In view of (i),
$R$ is finite. Since $\ang{x}$ is supercomplemented in $G$, for some
$T\leq G$, $G=RT$ and $F\cap T=1$. Since $|G:T|<\infty$, by
Poincare's Theorem, there exists $N\unlhd G$, $N\subseteq T$ such
that $|G:N|<\infty$. Obviously, $g\notin N$. In view of
arbitrariness of $g$, $G$ is residually finite.

(iv) Let $P$ be any finite $p$-subgroup of $G$. Since $F=\ang{P,x}$
is finite (see (i)), in view of Sylow's Theorem, in $F$ there exists
a Sylow $p$-subgroup $S$, containing $\ang{x}$ and some subgroup
$P^{*}$ conjugated with $P$. By Lemma 1, $\ang{x}$ is
supercomplemented in $S$.

If $p=2$, in view of Theorem \cite{Kreknin}, $d(S)\leq 3$. At the
same time, $d(P)\leq 3$. Consequently, each $2$-subgroup of $G$ is
solvable of derived length $\leq 3$.

Let $p\neq 2$. Then, in view of Theorem 2 \cite{Ch_Tr1}, $d(S)\leq
2$. Thus $d(P)\leq 2$ and each $p$-subgroup of $G$ is metabelian.

Further, in view of Theorem 1 \cite{Ch_Tr1}, $S$ has a normal
elementary abelian subgroup of index $\leq m!$. Then, $P^{*}$ and,
consequently, $P$ have a normal elementary abelian subgroup of index
$\leq m!$. Since $G$ is locally finite, with regard to Lemma 2, each
its $p$-subgroup is almost elementary abelian and nilpotent. Theorem
is proven.

\vspace{4mm}{\bfseries\itshape Proof of Theorem 2}. In view of
Proposition 4, $G$ is an $RN$-group. So by Theorem 1, the present
theorem is correct.

\vspace{4mm}{\bfseries\itshape Proof of Theorem 3}. Take $x\in
G\backslash H$ and $q\in \mathbb{P}$ such that $q \nmid
|\ang{x}:\ang{x}\cap H|$. Then $\ang{x^{q}}\nsubseteq H$. So
$\ang{x^{q}}$ is complemented in $G$. At the same time,
$\ang{x^{q}}$ is complemented in $\ang{x}$ (S.N.Chernikov's Lemma).
Hence, clearly, follows that $\ang{x}$ is of finite order.
Obviously, for some $p\in \mathbb{P}$, there is a $p$-element
$g\in\ang{x}\setminus H$. Fix $g$. Put $|\ang{g}|=m$. By Lemma 1,
$\ang{g}$ is supercomplemented in $G$.

Let $G_{\iota}$, $\iota\in I$, be all finitely generated subgroups
of $G$ containing $\ang{g}$, and let $G_{\iota}^{\varphi}$ be a
finite homomorphic image of some $G_{\iota}$.

If $G_{\iota}^{\varphi}\neq (H\cap G_{\iota})^{\varphi}$, then, in
consequence of Lemma 1, every subgroup of $G_{\iota}^{\varphi}$ not
belonging to $(H\cap G_{\iota})^{\varphi}$ is complemented in
$G_{\iota}^{\varphi}$. Therefore by Theorem \cite{Ch_Tr2},
$G_{\iota}^{\varphi}$ is solvable.

Let $G_{\iota}^{\varphi}=(H\cap G_{\iota})^{\varphi}$. Then
$\ang{g}^{\varphi}=1$. In view of Lemma 1, $\ang{g}^{\varphi}$ is
supercomplemented in $G_{\iota}^{\varphi}$. So all subgroups of
$G_{\iota}^{\varphi}$ are complemented in $G_{\iota}^{\varphi}$.
Therefore, in view of N.V.Chernikova's Theorem \cite{NVCher,NCher},
$G_{\iota}^{\varphi}$ is solvable.

Thus all finite $G_{\iota}^{\varphi}$ are solvable.

Further, since $\ang{g}^{\varphi}$ is a supercomplemented subgroup
of $G_{\iota}^{\varphi}$ (see Lemma 1), in view of Proposition 2,
$d(G_{\iota}^{\varphi})\leq 2n+3$ where $n=\max\{m(m-1)+m\log_{2}m,
1\}$. Then, in consequence of Lemma 3 \cite{Ch_Tr1}, $G_{\iota}$ is
finite and $d(G_{\iota})\leq 2n+3$. Therefore, obviously, $G$ is
solvable locally finite (and $d(G)\leq 2n+3$).

Then, in view of Theorem 1 (iii), $G$ is residually finite.

If $G$ is primary or all its primary subgroups are elementary
abelian, then, with regard to Theorem 1, the statement (iii) of the
present theorem is valid.

Let $G$ be non-primary and let for some $p\in \mathbb{P}$, $G$
contains some non-elementary abelian $p$-subgroup. Then, because of
$G$ is locally finite, it, obviously, contains some finite
non-elementary abelian $p$-subgroup $P$.

Take any $q\in \pi(G)\backslash\{p\}$, any finite $q$-subgroup
$Q\neq 1$ of $G$. Put $F=\ang{P,Q,g}$. Then, by proven above, $F$ is
finite solvable. Therefore for some $K\unlhd F$, ${\mathop {
\bigcap} \limits_{y\in F}}(F\cap H)^{y}\subseteq K$ and $|F:K|=r\in
\mathbb{P}$.

Take any $s\in \mathbb{P}\backslash\{r\}$ such that $s\mid |F|$ and
a Sylow $s$-subgroup $S$ of $F$. Then $1\neq S\subseteq K$.
Therefore, by Frattini argument (see, for instance, \cite{Rob}),
$F=KN_{F}(S)$. Take any Sylow $r$-subgroup $R$ of $N_{F}(S)$. Then
$R\nsubseteq K$. Since $R\nsubseteq {\mathop { \bigcap}
\limits_{y\in F}}(F\cap H)^{y}$, for some $u\in F$, $R\nsubseteq
(F\cap H)^{u}$. At the same time, $R^{u^{-1}}\nsubseteq H$.
Consequently, $R^{u^{-1}}$ is supercomplemented in $G$. Therefore
$R$ is supercomplemented in $G$. Then, in view of Lemma 1, $R$ is
supercomplemented in $RS$.

Let $T=\{b\in Z(S):b^{s}=1\}$ and let $L$ be a complement of $RT$ in
$RS$. Since $S$ is a normal Sylow $s$-subgroup of $RS$, obviously,
$L$ complements $T$ in $S$. Then, because of $T\subseteq Z(S)$,
clearly, $S=T\times L$. Therefore, inasmuch as $S$ is a finite
$s$-group and $L$ has no central in $S$ elements of order $s$,
$L=1$. So $T=S$, i.e. $S$ is elementary abelian. Since $F$ contains
the non-elementary abelian $p$-subgroup $P$ and for each $s\in
\mathbb{P}\backslash \{r\}$, a Sylow $s$-subgroup of $F$ is
elementary abelian, we have $p=r$. Thus all finite $q$-subgroups $Q$
of $G$ are elementary abelian. Therefore, obviously, the assertion
of (iii) of the present theorem concerning a $q$-subgroup is valid.

Take any $p$-element $a\in F\backslash K$. Then $\ang{a}$ is
supercomplemented in $G$. So, in consequence of Theorem 1, the
assertion of (iii) concerning a $p$-subgroup is valid. Theorem is
proven.

\vspace{4mm}{\bfseries\itshape Proof of Theorem 4}. In view of
Theorem 3, $\mathfrak{X}\subseteq \mathfrak{Y}$.

Clearly, $G=\ang{x}\leftthreetimes (\ang{a}\times \ang{b})$ with
$|\ang{x}|=8$, $|\ang{a}|=|\ang{b}|=2$ and $x^{a}=x^{-1}$,
$x^{b}=x^{5}$. Obviously, $\ang{x}$ is supercomplemented in $G$.
Indeed, if $\ang{x}\subseteq K\leq G$, then a complement $D$ to
$K\cap (\ang{a}\times \ang{b})$ in $\ang{a}\times \ang{b}$ is a
complement to $K$ in $G$. Thus $G\in \mathfrak{Y}$.

Show that $G$ has no $C$-separating subgroups, i.e. that $G\in
\mathfrak{Y}\backslash \mathfrak{X}$. Clearly, it is sufficient to
show that $G$ has no $C$-separating subgroups of index 2. Suppose
contrary.

Let $H$ be a $C$-separating subgroup of index 2. Obviously, each
subgroup of index 2 of $G$ contains the derived group
$G'=\ang{x^{2}}$ (of order 4). The quotient group $G/\ang{x^2}$ is
elementary abelian of order 8. Therefore it contains exactly 7
subgroups of order 4. Consequently, $G$ contains exactly 7 subgroups
of index 2. It is easy to see that the subgroups
$\ang{x}\leftthreetimes \ang{a}$, $\ang{x}\leftthreetimes \ang{b}$,
$\ang{x^2}\leftthreetimes (\ang{a}\times\ang{b})$,
$(\ang{x^2}\leftthreetimes \ang{b})\leftthreetimes\ang{xa}$,
$\ang{x}\leftthreetimes \ang{ab}$, $\ang{xb}\leftthreetimes
\ang{a}$, $(\ang{x^2}\leftthreetimes
\ang{ab})\leftthreetimes\ang{xa}$ are 7 pairwise distinct subgroups
of index 2 of $G$ complemented in it respectively by the subgroups
$\ang{b}$, $\ang{a}$, $\ang{xa}$, $\ang{x^2a}$, $\ang{a}$,
$\ang{b}$, $\ang{b}$. Thus $H$ is complemented in $G$. Then, in view
of (\cite{Spiv}, Theorem 2.1), $G=B\leftthreetimes A$ where $B$ is
elementary abelian and $|A|=2$. So an exponent of $G$ does not
exceed $4$, which is a contradiction (recall that $|\ang{x}|=8$).
Theorem is proven.

{\small

\renewcommand{\refname}{References}


}

\newpage

\end{document}